\documentclass[10pt]{amsart}
\usepackage{amssymb,latexsym,amsmath,amsfonts}
\usepackage[mathscr]{eucal}

\numberwithin{equation}{section}
\theoremstyle{definition}
\newtheorem{definition}{Definition}[section]
\newtheorem{remark}[definition]{Remark}
\newtheorem{example}[definition]{Example}
\newtheorem{custom}[definition]{}
\theoremstyle{plain}
\newtheorem{theorem}[definition]{Theorem}
\newtheorem{lemma}[definition]{Lemma}

\newtheorem{result}[definition]{Result}

\newcommand{\zetbar}{\overline{\zeta}}
\newcommand{\zt}{\zeta}

\newcommand{\Lmb}{\Lambda}
\newcommand{\lam}{\lambda}

\newcommand{\fee}{\varphi}

\newcommand{\rt}{\mathfrak{r}}
\newcommand{\banal}{\mathcal{A}}

\newcommand{\OM}{\Omega}
\newcommand{\Dsc}{\overline{D}}

\newcommand{\hol}{\mathcal{O}}

\newcommand{\poinc}{p_D}
\newcommand\Lem[1]{\widetilde{\kappa}_{#1}}

\newcommand\lmink[1]{{\rm h}_{\Lambda,{#1}}}
\newcommand{\seth}{\mathcal{S}}   
\newcommand\mink[1]{\mathcal{M}_{#1}}

\newcommand{\Cn}{\mathbb{C}^n}
\newcommand{\cplx}{\mathbb{C}}

\begin{document}

\title[Nonisotropically balanced domains and the Nevanlinna-Pick problem]{Nonisotropically balanced 
domains, Lempert \\
function estimates, and the spectral\\
Nevanlinna-Pick problem}
\author{Gautam Bharali}
\address{Department of Mathematics, Indian Institute of Science, Bangalore -- 560 012}
\email{bharali@math.iisc.ernet.in}

\begin{abstract} 
We introduce the notion of a {\em $\Lmb$-nonisotropically balanced domain} and show 
that the symmetrized polydisc in $\Cn, \ n\geq 2$, is an example of such a domain. Given 
a $\Lmb$-nonisotropically balanced domain $\OM$, we derive effective estimates from above
and from below for the Lempert function $\Lem{\OM}$ at $(0,z)\in\OM\times\OM$. 
We use these estimates to derive certain conditions for realising a two-point 
Nevanlinna-Pick interpolation in the symmetrized polydisc. Applying the ideas used in the
derivation of our Lempert function estimates to the so-called spectral unit ball $\OM_n$,
we deduce: {\em a)} a formula for the Lempert function at $(0,W)\in\OM_n\times\OM_n$; and
{\em b)} a necessary and sufficient condition for realising a two-point Nevanlinna-Pick
interpolation in the spectral unit ball.
\end{abstract}
\maketitle

\section{Introduction and statement of results}\label{S:intro}

This paper is partly motivated by the desire to obtain effective estimates for the Lempert
function $\Lem{\OM}$ for a domain $\OM\subset\Cn, \ n\geq 2$ (refer to Definition \ref{D:Lempert}
below). While effective formulae for $\Lem{\OM}$ are known for special classes of domains --- such
as when $\OM$ is a balanced domain or a Reinhardt domain --- estimates for $\Lem{\OM}$ are not known
even for many interesting examples in $\Cn$. It is with this situation in mind that we introduce a
new notion: that of a $\Lmb$-nonisotropically balanced domain. Our decision to single out this class 
of domains stems from this paper's slant towards the spectral Nevanlinna-Pick problem. These assertions
will be clearer once we have presented the following definition and an example.
\smallskip

\begin{definition} Let $\Lmb=(\lam_1,\dots,\lam_n)$ be an $n$-tuple of positive integers
that are relatively prime. A domain $\OM\subset\Cn, \ n\geq 2$, is said to be 
{\em $\Lmb$-nonisotropically balanced} if, whenever $z=(z_1,\dots,z_n)\in\OM$, then
$(\zt^{\lam_1}z_1,\dots,\zt^{\lam_n}z_n)\in\OM \ \forall\zt\in\Dsc$.
\end{definition}
\smallskip

In the above definition, and for the remainder of this paper, $D$ shall denote the open
unit disc in $\cplx$. We remark that in the terminology of the above definition, 
{\em balanced domains} in $\Cn$ are simply $(1,1,\dots,1)$-nonisotropically balanced
domains.
\smallskip

\begin{example} {\em The symmetrized polydisc in $\Cn$}

\noindent{The symmetrized polydisc in $\Cn$, denoted by $G_n$, is defined by
\[
G_n := 
\{(s_1,\dots,s_n)\in\Cn : \text{all the roots of 
$z^n-s_1z^{n-1}+\dots+(-1)^ns_n = 0$ lie in $D$}\}.
\]
$G_n$ is $(1,2,\dots,n)$-nonisotropically balanced. This follows from the fact that if
$\zt\in\Dsc$ and $\{r_1,\dots,r_n\}\subset D$ are the roots, repeated according to 
multiplicity, of
\[
z^n-s_1z^{n-1}+\dots+(-1)^ns_n \ = \ 0, \quad (s_1,\dots,s_n)\in G_n,
\]
then $\zt r_1,\dots,\zt r_n$ are the roots of
\[
z^n-(\zt s_1)z^{n-1}+\dots+(-1)^{n-1}(\zt^{n-1}s_{n-1})z
+(-1)^n(\zt^n s_n) \ = \ 0,
\]
and they all lie in $D$. In other words: 
$(s_1,\dots,s_n)\in G_n \Longrightarrow (\zt s_1,\zt^2s_2,\dots,\zt^ns_n)\in G_n$.
Hence $G_n$ is $(1,2,\dots,n)$-nonisotropically balanced.}
\end{example}
\smallskip

\begin{remark} It was pointed out to the author that, unbeknownst to him, the above definition
has appeared earlier in the preprint \cite{nikolov:spcedbcd05} by Nikolov. The terms
$\Lmb$-nonisotropically balanced domain, used herein, and $(k_1,\dots,k_n)$-balanced domain in
\cite{nikolov:spcedbcd05} are the same. Consequently, Lemma \ref{L:MinkPluri} below and
Prop.~1 of \cite{nikolov:spcedbcd05} are the same. The idea behind the argument presented in
both results seems to go back to Globevnik \cite{globevnik:Slsr74}. The ``only if'' part of
Theorem~\ref{T:specUB} below is a special case of Globevnik's result.
\end{remark}
\smallskip

The symmetrized polydisc has drawn quite a lot of attention lately owing to its 
connection with the spectral Nevanlinna-Pick problem. This problem is stated
as follows:
\smallskip

\begin{itemize}
\item[(*)] {\em Given $m$ distinct points $\zt_1,\dots,\zt_m\in D$ and matrices 
$W_1,\dots, W_m$ in the spectral unit ball $\OM_n:=\{W\in M_n(\cplx):r(W)<1\}$,
find conditions on $\{\zt_1,\dots,\zt_m\}$ and $\{W_1,\dots,W_m\}$ such that there 
exists a holomorphic map $F:D\longrightarrow\OM_n$ satisfying $F(\zt_j)=W_j, \
j=1,\dots,m$.}
\end{itemize}

\noindent{In the above statement, $r(W)$ denotes the spectral radius of the $n\times n$ 
matrix $W$. The papers \cite{aglerYoung:2psNPp00}, \cite{aglerYoung:2b2sNPp04} and
\cite{costara:22sNPp05} are just some of the recent papers dealing with the above
problem. Note that if $W$ is an $n\times n$ complex matrix, then $W\in\OM_n$ is 
equivalent to the fact that the coefficients of its characteristic polynomial determine
the coordinates of a point in the symmetrized polydisc. This observation forms the basis 
of recent investigations into the problem (*). This motivates another interesting 
interpolation problem analogous to (*), namely:}
\smallskip

\begin{itemize}
\item[(**)] {\em Given $m$ distinct points $\zt_1,\dots,\zt_m\in D$ and points
$p_1,\dots, p_m$ in the symmetrized polydisc $G_n$, find conditions on
$\{\zt_1,\dots,\zt_m\}$ and $\{p_1,\dots,p_m\}$ such that there exists a
holomorphic map $F:D\longrightarrow G_n$ satisfying $F(\zt_j)=p_j, \
j=1,\dots,m$.}
\end{itemize}

\noindent{One of the objectives of this paper is to show how estimates for the Lempert 
function $\Lem{\OM}$ can be used to derive:
\begin{itemize}
\item A necessary condition and a sufficient condition for the solvability of a $2$-point
interpolation problem in $G_n$ --- i.e., the interpolation problem (**) with $m=2$; and
\item A necessary and sufficient condition for the solvability of a $2$-point interpolation
problem in the spectral unit ball --- i.e., the interpolation problem (*) with $m=2$.
\end{itemize}}
\smallskip

Before we discuss these results, let us return to the basic issue of estimating
$\Lem{\OM}$ for a more general class of domains. We begin with the definition of the
Lempert function.
\smallskip

\begin{definition}[from \cite{lempert:mKrdb81}]\label{D:Lempert} 
Let $\OM$ be a domain in $\Cn, \ n\geq 2$, and let $z_1,z_2\in\OM$. The {\em Lempert function}
$\Lem{\OM}(z_1,z_2)$ is defined as
\begin{equation}\label{E:Lempert}
\Lem{\OM}(z_1,z_2) :=
\inf\{\poinc(0,\zt) : \zt\in D \ \text{and} \ \exists\fee\in\hol(D;\OM) \ \text{such
that} \ \fee(0)=z_1, \ \fee(\zt)=z_2 \},
\end{equation}
where $\poinc$ denotes the Poincar{\'e} distance on the unit disc.
\end{definition}
\smallskip

A comment on notation: given complex domains $X$ and $Y$, $\hol(X;Y)$ denotes the class of all
holomorphic mappings from $X$ into $Y$. The proof of our first result, which provides estimates
for $\Lem{\OM}(0,z)$, $z\in\OM$ and $\OM$ a $\Lmb$-nonisotropically balanced domain, exploits 
some of the ideas in the literature used in expressing $\Lem{G}(0,z)$ --- $G$ here being
a {\em balanced} domain ---  in terms of the Minkowski functional of $G$. For this, we would 
need a substitute for the Minkowski functional. We thus propose the following
\smallskip

\begin{definition} Let $\OM$ be a $\Lmb$-nonisotropically balanced domain in $\Cn$ and let 
$z\in\Cn$. The {\em Minkowski $\Lmb$-functional of $\OM$}, denoted $\lmink{\OM}$, is 
defined by
\[
\lmink{\OM}(z) := \inf\left\{t>0:\left(\frac{z_1}{t^{\lam_1}},\dots,
\frac{z_n}{t^{\lam_n}}\right)\in\OM\right\}.
\]
\end{definition}
\smallskip

We now have all the elements necessary to state our first result.
\smallskip

\begin{theorem}\label{T:LemMink}
Let $\OM$ be a $\Lmb$-nonisotropically balanced pseudoconvex domain, and let $\lmink{\OM}$ be
its Minkowski $\Lmb$-functional. Let $L:=\max\{\lam_j:j=1,\dots,n\}$. Then, for any $z\in\OM$,
the Lempert function $\Lem{\OM}(z)$ satisfies
\begin{equation}\label{E:LemEst}
\tanh^{-1}(\lmink{\OM}(z)^L) \ \leq \ \Lem{\OM}(0,z) \ \leq \ \tanh^{-1}(\lmink{\OM}(z)).
\end{equation}
\end{theorem}
\smallskip

Note that when $\OM$ is a balanced domain, then 
$\Lem{\OM}(0,z)=\tanh^{-1}(\lmink{\OM}(z))
=\tanh^{-1}(\mink{\OM}(z))$ $\forall z\in\OM$, where $\mink{\OM}$ is the
Minkowski functional of $\OM$. This equation has, of course, long been an established fact. We
borrow some of the ideas used in proving this equality to establish Theorem \ref{T:LemMink}. 
This will require studying the properties of
the Minkowski $\Lmb$-functional of a nonisotropically balanced domain. These properties are
investigated in Section \ref{S:MinkProps}. The proof of Theorem \ref{T:LemMink} is given at the
end of that section.
\smallskip

Note that simply from the definition \eqref{E:Lempert}, we have the following Schwarz lemma for
the domain $\OM$:
\begin{equation}\label{E:Schwarz}
\left.\begin{array}{rl}
f &\in \ \hol(D;\OM),\\
f(\zt_j) &= \ z_j, \ j=1,2,\\
\zt_j &\in \ D,\\
z_j &\in \ \OM\end{array}\right\} \ \Longrightarrow \ 
\Lem{\OM}(z_1,z_2)  \ \leq \ \poinc(\zt_1,\zt_2).
\end{equation}
Unfortunately, the above statement conveys very little information unless one knows 
$\Lem{\OM}$ explicitly. The Lempert function of the symmetrized bidisc $G_2$ is exactly known, 
whence one has an {\em explicit} Schwarz lemma for the symmetrized bidisc. The reader is 
referred to \cite[Theorem 1.1]{aglerYoung:Slfsb01} by Agler \& Young. In higher dimensions,
a necessary condition for the type of interpolation described in the
hypothesis of \eqref{E:Schwarz}, with $\OM=G_n, \ n\geq 3$, has been established in the
recent paper \cite{costara:osNPp05} by Costara. However, to the best of our knowledge,
\begin{itemize}
\item Estimates from both above and below are yet unknown for the Lempert function of
$G_n, \ n\geq 3$; and
\item No sufficient conditions are known for the solvability of the two-point interpolation
problem given in the hypothesis of \eqref{E:Schwarz} when $\OM=G_n, \ n\geq 3$.
\end{itemize}
To this end, we provide a necessary condition and a sufficient condition for the solvability
of the interpolation problem under discussion.
\smallskip

\begin{theorem}\label{T:2pointSym}
Given any point $s:=(s_1,\dots,s_n)$ in the symmetrized polydisc $G_n$, define the
rational function
\begin{equation}\label{E:ratImp}
F_s(z) \ := \ \frac{(-1)^n ns_n z^{n-1}+(-1)^{n-1}(n-1)s_{n-1}z^{n-2}+\dots+(-s_1)}
		{n-(n-1)s_1 z+\dots+(-1)^{n-1}s_{n-1}z^{n-1}}, \quad z\in\cplx.
\end{equation}
Let $\zt_1,\zt_2\in D$ and $s\in G_n$. Then:
\begin{enumerate}
\item[1)] If there exists a map $f\in\hol(D;G_n)$ such that $f(\zt_1)=0$ and $f(\zt_2)=s$, then
\begin{equation}\label{E:2pointNece}
\max\left\{\inf
\left\{t^n: t>0 \ \text{and} \ \frac{1}{t}\sup_{|z|=1/t}|F_s(z)|<1\right\}, \
\sup_{|z|=1}|F_s(z)| \right\} \ \leq
\ \left|\frac{\zt_2-\zt_1}{1-\zetbar_1\zt_2}\right|.
\end{equation}
\item[2)] If
\begin{equation}\label{E:2pointSuff}
\inf\left\{t>0:\frac{1}{t}\sup_{|z|=1/t}|F_s(z)|<1\right\} \ \leq
\ \left|\frac{\zt_2-\zt_1}{1-\zetbar_1\zt_2}\right|,
\end{equation}
then there exists a map $f\in\hol(D;G_n)$ such that $f(\zt_1)=0$ and $f(\zt_2)=s$.
\end{enumerate}
\end{theorem}

Now consider the spectral unit ball, which is a balanced domain. If one could show that 
$\OM_n$ is pseudoconvex, then one would have an exact expression for the Lempert function
at $(0,W)\in\OM_n\times\OM_n$. One is able to show pseudoconvexity using Vesentini's theorem
\cite{vesentini:spr68}. Consequently, one obtains an analogue of Theorem \ref{T:2pointSym}
for the spectral unit ball. The approaches to solving the problem (*) that are discussed in 
the aforementioned papers depend on using information about the problem (**) to analyse (*).
Since the coefficients of the characteristic polynomial of any $W\in M_n(\cplx)$ do not
alone encode all the information about the Jordan structure of $W$, these approaches have
tackled (*) under the restriction that $W_1,\dots,W_n$ be non-derogatory (i.e., each $W_j$
is similar to its companion matrix). In contrast --- even though we address only a special
case of (*) with $m=2$ and $W_1=0$ --- our technique of calculating the Lempert function 
imposes no restrictions on the Jordan structure of $W_2\in\OM_n$. Before stating the 
pertinent result, let us fix the following notation:
\[
\Lem{n\times n}(W_1,W_2)  :=  \text{the Lempert function evaluated at 
	$(W_1,W_2)\in\OM_n\times\OM_n$}.
\]
We can now state our final theorem.
\smallskip

\begin{theorem}\label{T:specUB}
For any matrix $W\in\OM_n$, let 
\[
z^n-s_1z^{n-1}+\dots+(-1)^{n-1}s_{n-1}z+(-1)^n s_n \ = \ 0
\]
denote the characteristic equation of $W$, and define the rational function
\[
F_{s(W)}(z) \ := \ \frac{(-1)^n ns_n z^{n-1}+(-1)^{n-1}(n-1)s_{n-1}z^{n-2}+\dots+(-s_1)}
                {n-(n-1)s_1 z+\dots+(-1)^{n-1}s_{n-1}z^{n-1}}, \quad z\in\cplx.
\]
\begin{enumerate}
\item[1)] Let $W\in\OM_n$. Then,
\[
\Lem{n\times n}(0,W) \ = \
\tanh^{-1}\left[ \inf\left\{t>0:\frac{1}{t}\sup_{|z|=1/t}|F_{s(W)}(z)|<1\right\}\right] \ 
= \ \tanh^{-1}(r(W)).
\]
\item[2)] Let $\zt_1,\zt_2\in D$ and $W\in\OM_n$.
There exists a map $f\in\hol(D;\OM_n)$ such that $f(\zt_1)=0$ and $f(\zt_2)=W$
if and only if
\begin{equation}\label{E:2pointUBiff}
r(W) \ = \ \inf\left\{t>0:\frac{1}{t}\sup_{|z|=1/t}|F_{s(W)}(z)|<1\right\} \ \leq
\ \left|\frac{\zt_2-\zt_1}{1-\zetbar_1\zt_2}\right|.
\end{equation}
\end{enumerate}
\end{theorem}
\smallskip             

\section{Properties of $\lmink{\OM}$}\label{S:MinkProps}

We begin with the following elementary result:
\smallskip

\begin{lemma}\label{L:MinkProps}
 Let $\OM$ be a $\Lmb$-nonisotropically balanced domain in $\Cn, \ n\geq 2$, and
$\OM\neq\Cn$. Then:
\begin{enumerate}
\item[a)] $\lmink{\OM}(\zt^{\lam_1}z_1,\dots,\zt^{\lam_n}z_n) \ = \ 
|\zt|\lmink{\OM}(z) \ \forall\zt\in\cplx$ and $\forall z\in\Cn$.
\item[b)] $\OM=\{z\in\Cn:\lmink{\OM}(z)<1\}$.
\end{enumerate}
\end{lemma}
\begin{proof}
Define the set
\[
\seth_{\Lmb}(z) \ := \ \left\{t>0:\left(\frac{z_1}{t^{\lam_1}},\dots,\frac{z_n}{t^{\lam_n}}\right)
\in\OM\right\}.
\]
Clearly, since $0$ is an interior point of $\OM$, $\seth_{\Lmb}(z)\neq\emptyset$. If 
$t_0\in\seth_{\Lmb}(z)$ and $s>t_0$, then
\[
\left(\frac{z_1}{s^{\lam_1}},\dots,\frac{z_n}{s^{\lam_n}}\right) \ = \ 
\left(\left(\frac{t_0}{s}\right)^{\lam_1}\frac{z_1}{t_0^{\lam_1}},\dots,
\left(\frac{t_0}{s}\right)^{\lam_n}\frac{z_n}{t_0^{\lam_n}}\right)\in\OM
\]
because $\OM$ is $\Lmb$-nonisotropically balanced. This implies that
\begin{equation}\label{E:t-range/h}
\seth_{\Lmb}(z) \ = \ (\lmink{\OM}(z),+\infty).
\end{equation}
Now note that for $\zt\in\cplx$,
\begin{align}
{} & \ t\in\seth_{\Lmb}(\zt^{\lam_1}z_1,\dots,\zt^{\lam_n}z_n) \notag\\
\Leftrightarrow & \ \left(\frac{\zt^{\lam_1}z_1}{t^{\lam_1}},\dots,
	\frac{\zt^{\lam_n}z_n}{t^{\lam_n}}\right)\in\OM \notag\\
\Leftrightarrow & \ \left(e^{-i{\lam_1}{\rm Arg}(\zt)}\frac{\zt^{\lam_1}z_1}{t^{\lam_1}},\dots,
        e^{-i{\lam_n}{\rm Arg}(\zt)}\frac{\zt^{\lam_n}z_n}{t^{\lam_n}}\right)\in\OM 
	\notag\\
\Leftrightarrow & \  \left(\frac{z_1}{(t/|\zt|)^{\lam_1}},\dots,
        \frac{z_n}{(t/|\zt|)^{\lam_n}}\right)\in\OM \notag\\
\Leftrightarrow & \ t/|\zt|\in\seth_{\Lmb}(z)\notag
\end{align}
This implies that $\seth_{\Lmb}(\zt^{\lam_1}z_1,\dots,\zt^{\lam_n}z_n)=(|\zt|\lmink{\OM}(z),+\infty)$.
In view of \eqref{E:t-range/h}, we conclude that 
$\lmink{\OM}(\zt^{\lam_1}z_1,\dots,\zt^{\lam_n}z_n) \ = \
|\zt|\lmink{\OM}(z)$.
\smallskip

To prove (b), first note that if $z\in\Cn$ and $\lmink{\OM}(z)<1$, then $1\in\seth_{\Lmb}(z)$, i.e.,
$z\in\OM$. This means that $\{z\in\Cn:\lmink{\OM}(z)<1\}\subseteq\OM$. To establish the opposite 
inclusion, we first consider

\noindent{{\bf {\em Case 1}.} {\em $z\in\OM$ and $\lmink{\OM}(z)=T>1$}}

\noindent{In this case, there exists an $r\in(1,T)$ such that $r\notin\seth_{\Lmb}(z)$, or 
equivalently
\[
Z \ := \ \left(\frac{z_1}{r^{\lam_1}},\dots,\frac{z_n}{r^{\lam_n}}\right)\notin\OM.
\]
But, as $(1/r)<1$ and $z\in\OM$, $Z$ must belong to $\OM$, which is a contradiction.}
\smallskip

\noindent{{\bf {\em Case 2}.} {\em $z\in\OM$ and $\lmink{\OM}(z)=1$}}

\noindent{In this case, using the argument in Case 1, we can infer that 
$z\in\overline{\Cn\setminus\OM}$. This implies that $z\notin\OM$, which is 
again a contradiction.}
\smallskip

We have just shown that $\OM\bigcap\{z\in\Cn:\lmink{\OM}\geq 1\}=\emptyset$. In conjunction
with the earlier inclusion, we get  $\OM=\{z\in\Cn:\lmink{\OM}(z)<1\}$. 
\end{proof}
\smallskip

Our next lemma establishes a crucial fact about $\lmink{\OM}$: namely that it is a
plurisubharmonic function on $\Cn$ if $\OM$ is pseudoconvex. Once this is established,
we can to exploit ideas that have been used in the study of invariant metrics on
balanced domains. To prove this lemma, we will require the following result (refer,
for instance, to Appendix PSC in \cite{jarnickiPflug:idmca93} by Jarnicki \& Pflug):
\smallskip

\begin{result}\label{R:plurish} 
Let $G$ be a balanced domain in $\Cn, \ n\geq 2$, and let
\[
\mink{G}(z) := \inf\left\{t>0:\left(\frac{z_1}{t},\dots,
\frac{z_n}{t}\right)\in G\right\}, \quad z\in\Cn,
\]
be the Minkowski functional of $G$. $G$ is pseudoconvex if and only if $\mink{G}$
is plurisubharmonic on $\Cn$.
\end{result}
\smallskip

Using this result, we can now prove the following
\smallskip

\begin{lemma}\label{L:MinkPluri}
Let $\OM$ be a $\Lmb$-nonisotropically balanced domain in $\Cn, \ n\neq 2$, and $\OM\neq\Cn$. If 
$\OM$ is pseudoconvex, then $\lmink{\OM}$ is plurisubharmonic on $\Cn$.
\end{lemma}
\begin{proof}
It is quite evident the $\lmink{\OM}$ is upper semicontinuous (for instance, $\lmink{\OM}$ can
be rewritten as the lower envelope of a $1$-parameter family of upper semicontinuous functions).
Now, let us define two auxilliary objects:
\begin{align}
\omega \ &:= \ \{z\in\Cn:(z_1^{\lam_1},\dots,z_n^{\lam_n})\in\OM\},\label{E:auxdom} \\
{\rm h}^{(\OM)}(z) \ &:= \ \inf\left\{t>0:\left(\frac{z_1^{\lam_1}}{t^{\lam_1}},\dots,
\frac{z_n^{\lam_n}}{t^{\lam_n}}\right)\in \OM\right\}, \quad z\in\Cn. \notag
\end{align}
Observe that $\omega$ is a balanced domain in $\Cn$. Note, furthermore, that
\begin{align}
\mink{\omega}(z) \ &= \ \inf\left\{t>0:\left(\frac{z_1}{t},\dots,
\frac{z_n}{t}\right)\in \omega\right\} && {} \notag\\
&= \ \inf\left\{t>0:\left(\frac{z_1^{\lam_1}}{t^{\lam_1}},\dots,
\frac{z_n^{\lam_n}}{t^{\lam_n}}\right)\in \OM\right\} && [\text{follows from
definition \eqref{E:auxdom} above}] \notag \\
&= \ {\rm h}^{(\OM)}(z). && {}\label{E:auxeq}
\end{align}
Now consider the map
\[
P_{\Lmb} \ : \ (z_1,\dots,z_n)\longmapsto(z_1^{\lam_1},\dots,z_n^{\lam_n}).
\]
Observe that $P_{\Lmb}:\omega\longrightarrow\OM$ is a proper mapping onto $\OM$ and, in fact,
$P_{\Lmb}^{-1}(\OM)=\omega$. Since $P_{\Lmb}$ is a proper holomorphic $\omega\longrightarrow\OM$
mapping, and $\OM$ is pseudoconvex, $\omega$ is also pseudoconvex. 
\smallskip

Let us define the algebraic variety
\[
V \ := \ \bigcup_{j=1}^n\{z\in\Cn:z_j=0\},
\]
and the function 
\[
u(w) \ := \ \frac{1}{\lam_1\lam_2\dots\lam_n}\sum_{z\in P_{\Lmb}^{-1}\{w\}}{\rm h}^{(\OM)}(z) 
\quad\forall w\in\Cn\setminus V.
\]
Thus far, we have shown that
\begin{itemize}
\item $\omega$ is a balanced domain in $\Cn$;
\item $\omega$ is pseudoconvex; and
\item $\mink{\omega}(z)={\rm h}^{(\OM)}(z) \ \forall z\in\Cn$.
\end{itemize}
Combining these facts with Result \ref{R:plurish} we conclude that
\begin{equation}\label{E:keypluri}
u \ \text{\em is plurisubharmonic on $\Cn\setminus V$.}
\end{equation}
Note that for any $w\in\Cn\setminus V$, $W$ has $\lam_1\lam_2\dots\lam_n$ pre-images under
$P_{\Lmb}$. From this and the definition of ${\rm h}^{(\OM)}$, we can express $\lmink{\OM}$
in the following manner:
\[
\lmink{\OM}(w) \ = \ \begin{cases}
			u(w), & \text{if $w\in\Cn\setminus V$}, \\
			\lmink{\OM}(w), & \text{if $w\in V$}.
			\end{cases}
\]
In view of \eqref{E:keypluri}, and because
\begin{itemize}
\item $V$ is a closed pluripolar subset of $\Cn$; and
\item $\lmink{\OM}$ is locally bounded at each $w\in V$,
\end{itemize}
the removable singularities lemma for plurisubharmonic
functions implies that $\lmink{\OM}$ is plurisubharmonic on $\Cn$.
\end{proof}
\smallskip

\begin{custom}\label{Cu:proofLemMink}
\begin{proof}[{\bf The proof of Theorem \ref{T:LemMink}}]
We shall use $D^*$ to denote the punctured unit disc in $\cplx$. Recall that 
$L:=\max\{\lam_j:j=1,\dots,n\}$. For each $\zt\in D^*$, let 
$\rt_1(\zt),\dots,\rt_L(\zt)$ denote the distinct $L$th roots of $\zt$. Now
consider a $\fee\in\hol(D;\OM)$ such that $\fee(0)=0$ and such that there
exists a $\sigma\in D$ for which $\fee(\sigma)=z$. Since $\fee(0)=0$, we can
express $\fee$ as
\[
\fee(\zt) \ = \ (\zt\Phi_1(\zt),\dots,\zt\Phi_n(\zt))\quad\zt\in D,
\]
where $(\Phi_1,\dots,\Phi_n)\in\hol(D;\Cn)$. Notice that in view of 
Lemma \ref{L:MinkProps}, for any $\zt\in D^*$ and any $j=1,\dots,L$, we get
\begin{equation}\label{E:homo}
\lmink{\OM}\circ\fee(\zt) \ = \ |\zt|^{1/L}\lmink{\OM}
(\rt_j(\zt)^{L-\lam_1}\Phi_1(\zt),\dots,\rt_j(\zt)^{L-\lam_n}\Phi_n(\zt))
\quad(\zt\in D^*).
\end{equation}
Since we have shown that $\lmink{\OM}\in {\rm psh}(\Cn)$ (Lemma \ref{L:MinkPluri} above),
we can conclude, for any $\zt_0\in D^*$ and any sufficiently small neighbourhood $W(\zt_0)$
of $\zt_0$, that
\begin{equation}\label{E:keysubh}
W(\zt_0)\ni\zt\longmapsto
\lmink{\OM}(\rt_j(\zt)^{L-\lam_1}\Phi_1(\zt),\dots,\rt_j(\zt)^{L-\lam_n}\Phi_n(\zt)) \
\text{\em is subharmonic on $W(\zt_0)$.}
\end{equation}
Since subharmonicity is a local property, and since the function $U$ defined as
\[
U(\zt) \ := \ 
\sum_{j=1}^L\lmink{\OM}(\rt_j(\zt)^{L-\lam_1}\Phi_1(\zt),\dots,
\rt_j(\zt)^{L-\lam_n}\Phi_n(\zt)) \quad\forall\zt\in D^*
\]
is upper semicontinuous on $D^*$, $U$ is subharmonic on $D^*$. Finally, as $U$
is bounded in a punctured neighbourhood of $\zt=0$, $U$ extends to a subharmonic
function on $D$. From \eqref{E:homo} and Lemma \ref{L:MinkProps}/(b), we infer, for
each $r\in(0,1)$, that
\[
r^{1/L}U(\zt) \ = \ L\lmink{\OM}\circ\fee(\zt) \ < \ L\quad
\forall\zt:|\zt|=r.
\]
Therefore, by the Maximum Principle for subharmonic functions
\[
U(\zt) \ \leq \ L \qquad\forall\zt\in D.
\]
This tells us that
\[
L\lmink{\OM}(z) \ = \ L\lmink{\OM}\circ\fee(\sigma) \ = \ |\sigma|^{1/L}U(\sigma) \ 
\leq \ L|\sigma|^{1/L}.
\]
By Definition \ref{D:Lempert}, we see that 
\begin{equation}\label{E:1stHalf}
\Lem{\OM}(0,z) \ \geq \ \poinc(0,\lmink{\OM}(z)^L) \ = \ \tanh^{-1}(\lmink{\OM}(z)^L).
\end{equation}

To prove the other inequality in \eqref{E:LemEst}, we consider the following cases:
\smallskip

\noindent{{\bf {\em Case 1}.} {\em $z\in\OM$ is such that $\lmink{\OM}(z)\neq 0$.}}

\noindent{In this case, define
\[
\fee(\zt) \ := \ \left(\left(\frac{\zt}{\lmink{\OM}(z)}\right)^{\lam_1}z_1,\dots,
\left(\frac{\zt}{\lmink{\OM}(z)}\right)^{\lam_n}z_n\right).
\]
Note that
\[
\lmink{\OM}(z) \ > \ |\zt|\lmink{\OM}(z) \ = \ 
\lmink{\OM}(\zt^{\lam_1}z_1,\dots,\zt^{\lam_n}z_n)\quad\forall\zt\in D.
\]
Thus, by the definition of $\lmink{\OM}(\zt^{\lam_1}z_1,\dots,\zt^{\lam_n}z_n)$:
\[
 \left(\left(\frac{\zt}{\lmink{\OM}(z)}\right)^{\lam_1}z_1,\dots,
\left(\frac{\zt}{\lmink{\OM}(z)}\right)^{\lam_n}z_n\right)\in\OM\quad
\forall\zt\in D.
\]
Therefore, $\fee(D)\subset\OM$ and $\fee(\lmink{\OM}(z))=z$. Therefore, by
definition
\begin{equation}\label{E:2ndHalf1}
\Lem{\OM}(0,z) \ \leq \ \tanh^{-1}(\lmink{\OM}(z)).
\end{equation}}
\smallskip

\noindent{{\bf {\em Case 2}.} {\em $z\in\OM$ is such that $\lmink{\OM}(z)=0$.}}

\noindent{In this situation, for each $t\in(1,\infty)$, define
\[
\fee_t(\zt) \ := \ ((t\zt)^{\lam_1}z_1,\dots,(t\zt)^{\lam_n}z_n)
\]
Notice that, by assumption, $\lmink{\OM}(\fee_t(\zt))=|t\zt|\lmink{\OM}(z)=0$. By
Lemma \ref{L:MinkProps}/(b), $\fee_t(D)\subset\OM \ \forall t>1$. Furthermore
$\fee_t(1/t)=z$. Then,
\begin{equation}\label{E:2ndHalf2}
\Lem{\OM}(0,z) \ \leq \ \lim_{t\to\infty}\poinc(0,1/t) \ = \ 0 \ = \ 
\tanh^{-1}(\lmink{\OM}(z)).
\end{equation}}
\smallskip

From \eqref{E:1stHalf}, \eqref{E:2ndHalf1} and \eqref{E:2ndHalf2}, the result follows.
\end{proof}
\end{custom}

\section{The proofs of Theorems \ref{T:2pointSym} and \ref{T:specUB}}
\smallskip

To prove Theorem \ref{T:2pointSym}, we use the characterisation given below for a 
point $(s_1,\dots,s_n)$ to belong to the symmetrized polydisc. We point out that the
result below is {\em not the only characterisation} available for $s$ to belong to 
$G_n$; recall, for instance, the Schur-Cohn characterisation. However, since the
Schur-Cohn characterisation is no easier to check, in the context of Theorems
\ref{T:2pointSym} and \ref{T:specUB}, than the one presented below, we prefer to use 
the following criterion --- which has the advantage that it  could conceivably be used
to investigate the Carath{\'e}odory metric on $G_n$.
\smallskip

\begin{result}[Costara, \cite{costara:osNPp05}]\label{R:costaraChar}
Let $s=(s_1,\dots,s_n)\in\Cn$ and let $F_s(z)$ be the rational function given by
\eqref{E:ratImp}. Then, $s\in G_n$ if and only if $\sup_{|z|=1}|F_s(z)|<1$.
\end{result}

Using this result, we can now provide

\begin{custom}
\begin{proof}[{\bf The proof of Theorem \ref{T:2pointSym}}]
Assume that we are given $\zt_1,\zt_2\in D$ and $s\in G_n$, and that there exists a map 
$f\in\hol(D;G_n)$ such that $\zt_1=0$ and $\zt_2=s$. Since the biholomorphisms of $D$
act transitively on $D$, and since the Poincar{\'e} metric is invariant under 
biholomorphisms, by definition
\[
\Lem{\OM}(0,z) \ \leq \ \poinc(\zt_1,\zt_2) \ = \ \tanh^{-1}
			\left|\frac{\zt_2-\zt_1}{1-\zetbar_1\zt_2}\right|.
\]
Recall that $G_n$ is a $(1,2,\dots,n)$-balanced domain; in this proof, therefore,
$\Lmb$ will always denote the $n$-tuple $(1,2,\dots,n)$. Let us use the notation
$\lmink{n}$ to denote the Minkowski $(1,2,\dots,n)$-functional of $G_n$. In view of 
\eqref{E:LemEst}, the above implies that
\begin{equation}\label{E:intermed}
\tanh^{-1}(\lmink{n}(z)^n) \ \leq \ \tanh^{-1}
                        \left|\frac{\zt_2-\zt_1}{1-\zetbar_1\zt_2}\right|.
\end{equation}
Given $s\in G_n$ and $t>0$, let us define $t\bullet s:=(s_1/t,s_2/t^2,\dots,s_n/t^n)$.
We compute that
\begin{equation}\label{E:F_t*s}
F_{t\bullet s}(z) \ = \ \frac{1}{t}F_s\left(\frac{z}{t}\right).
\end{equation} 
Using Result \ref{R:costaraChar}, we get
\begin{align}
\lmink{n}(s) \ &= \ \inf\left\{t>0:\left(\frac{s_1}{t},\frac{s_2}{t^2}\dots,
\frac{s_n}{t^n}\right)\in G_n\right\} && {}\notag\\
&= \ \inf\left\{t>0:\sup_{|z|=1}|F_{t\bullet s}(z)|<1\right\} && 
[\text{from Result \ref{R:costaraChar}}] \notag\\
&= \ \inf\left\{t>0:\frac{1}{t}\sup_{|z|=1/t}|F_s(z)|<1\right\} && 
[\text{from \eqref{E:F_t*s}}].\notag
\end{align}
Using the above calculation in conjunction with \eqref{E:intermed}, we get
\begin{equation}\label{E:necess1}
\inf\left\{t^n:t>0 \ \text{and} \ \frac{1}{t}\sup_{|z|=1/t}|F_s(z)|<1\right\} \ 
\leq \ \left|\frac{\zt_2-\zt_1}{1-\zetbar_1\zt_2}\right|.
\end{equation}
As an easy corollary to Result \ref{R:costaraChar}, we also have the following
necessary condition for the existence of $f$ with the aforementioned properties
\cite[Corollary 3.1]{costara:osNPp05}:
\begin{equation}\label{E:necess2}
\sup_{|z|=1}|F_s(z)| \ \leq \ \left|\frac{\zt_2-\zt_1}{1-\zetbar_1\zt_2}\right|.
\end{equation}
From \eqref{E:necess1} and \eqref{E:necess2}, we conclude that
\begin{equation}
\max\left\{\inf
\left\{t^n: t>0 \ \text{and} \ \frac{1}{t}\sup_{|z|=1/t}|F_s(z)|<1\right\}, \
\sup_{|z|=1}|F_s(z)| \right\} \ \leq
\ \left|\frac{\zt_2-\zt_1}{1-\zetbar_1\zt_2}\right|.
\end{equation}   
\smallskip

Now assume \eqref{E:2pointSuff}. From the discussion in the preceding paragraph, we see
that this is the same as
\begin{equation}\label{E:implic}
\lmink{n}(s) \ \leq \ \left|\frac{\zt_2-\zt_1}{1-\zetbar_1\zt_2}\right|.
\end{equation}   
Since $G_n$ is a bounded domain, and $s\neq 0$, $\lmink{n}(s)>0$. Refer to the
argument given under Case 1 in the proof of Theorem \ref{T:LemMink}. That argument
establishes that the holomorphic map
\[
\fee(\zt) \ := \ \left(\left(\frac{\zt}{\lmink{n}(s)}\right)s_1,
\left(\frac{\zt}{\lmink{n}(s)}\right)^2s_2,\dots,
\left(\frac{\zt}{\lmink{n}(s)}\right)^ns_n\right)
\]
satisfies $\fee(D)\subset G_n$, $\fee(0)=0$ and $\fee(\lmink{n}(s))=s$. In view
of the inequality \eqref{E:implic}, the classical Nevanlinna-Pick theorem for the
unit disc in $\cplx$ tells us that there exists a $\psi\in\hol(D;D)$ such that 
$\psi(\zt_1)=0$ and $\psi(\zt_2)=\lmink{n}(s)$. Now, let us define $f:=\fee\circ\psi$. 
Clearly, $f$ has the properties asserted in Theorem \ref{T:2pointSym}/(2).
\end{proof}
\end{custom}
\smallskip

In view of Theorem \ref{T:LemMink}, the first part of Theorem \ref{T:specUB} would
follow quite simply if we could show that the spectral unit ball $\OM_n$ is a
pseudoconvex domain. To accomplish this, we would need the following result by
Vesentini:
\smallskip

\begin{result}[Vesentini, \cite{vesentini:spr68}]\label{R:vesentiniSpec} 
Let $\banal$ be a complex, unital Banach
algebra and let $r(x)$ denote the spectral radius of any element $x\in\banal$.
Let $f\in\hol(D;\banal)$. The the function $\zt\longmapsto r(f(\zt))$ is 
subharmonic on $D$.
\end{result}

Finally, we can provide
\smallskip

\begin{custom}
\begin{proof}[{\bf The proof of Theorem \ref{T:specUB}}]
In the notation of Result \ref{R:vesentiniSpec}, let $\banal=M_n(\cplx)$, which is
a unital Banach algebra with respect to the operations of matrix addition and matrix
multiplication. Consider any analytic disc $f\in\hol(D;\banal=M_n(\cplx))$. By 
Vesentini's theorem, $r\circ f$ is a subharmonic function --- where $r$ is the
function that maps each matrix in $M_n(\cplx)$ to its spectral radius. Since this
is true for any arbitrary analytic disc $f$, $r$ is, by definition, a 
plurisubharmonic function. Thus, $\OM_n$, which is defined as 
$\OM_n:=\{W\in M_n(\cplx):r(W)<1\}$, is a pseudoconvex domain. It is, of course,
a balanced domain. In this proof, therefore, $\Lmb$ will always denote the
$n$-tuple $(1,1,\dots,1)$. Hence, by Theorem \ref{T:LemMink} --- the following 
fact has also been established earlier
by slightly different methods --- we get:
\begin{equation}\label{E:minkLemEq}
\tanh^{-1}(\lmink{n\times n}(W)) \ = \ \Lem{n\times n}(0,W) \quad\forall W\in\OM_n,
\end{equation}
where $\Lmb=(1,1,\dots,1)$. It is evident that $h_{\Lambda,n\times n}(W)=r(W)$. Let us, however,
provide another computation for $h_{\Lambda,n\times n}(W)$. In the case when it is known that 
$W\in\Omega_n$, but its spectral radius is not explicitly known, this alternative expression 
is stabler to compute when $W$ is of large dimension.
For each $W\in M_n(\cplx)$, let $s_j(W), \ j=1,\dots,n$, be the coefficient of the
$z^{n-j}$-term of the characteristic polynomial of $W$. Referring to the notation in 
the statement of Theorem \ref{T:specUB}: let $s(W):=(s_1(W),\dots,s_n(W))$.
Observe that
\[
W\in\OM_n \ \Leftrightarrow \ s(W)\in G_n \ \Leftrightarrow \ \sup_{|z|=1}|F_{s(W)}(z)|<1.
\]
In view of the discussion in the proof of the previous theorem, the above implies:
\begin{align}
W/t\in\OM_n \ \Leftrightarrow \ t\bullet s(W)\in G_n \ &\Leftrightarrow \ 
\sup_{|z|=1}|F_{t\bullet s(W)}(z)|<1 \notag\\
&\Leftrightarrow \ \frac{1}{t}\sup_{|z|=1/t}|F_{s(W)}(z)|<1.\notag
\end{align}
Combining this with \eqref{E:minkLemEq}, we get
\begin{equation}\label{E:LempertSpecUB}
\Lem{n\times n}(0,W) \ = \ 
\tanh^{-1}\left[ \inf\left\{t>0:\frac{1}{t}\sup_{|z|=1/t}|F_{s(W)}(z)|<1\right\}\right].
\end{equation}
\smallskip

If there exists a map $f\in\hol(D;\OM_n)$ such that $f(\zt_1)=0$ and $f(\zt_2)=W$, then
simply by the definition
\[
\Lem{n\times n}(0,W) \ \leq \ \poinc(\zt_1,\zt_2).
\]
Thus, \eqref{E:2pointUBiff} follows from the formula \eqref{E:LempertSpecUB}. Now, conversely,
assume that for $\zt_1,\zt_2\in D$ and $W\in\OM_n$, the inequality \eqref{E:2pointUBiff} holds.
Since $\OM_n$ is {\em unbounded}, we will have to be careful. The existence of an $f\in\hol(D;\OM_n)$
will have to be analysed under the following two cases:
\smallskip

\noindent{{\bf {\em Case 1}.} {\em $W\in\OM_n$ is such that $\lmink{n\times n}(W)\neq 0$.}}

\noindent{The existence of an $f\in\hol(D;\OM_n)$ is completely analogous to the argument given
in the latter half of the proof of Theorem \ref{T:2pointSym}. We shall, therefore, not repeat
that argument.}
\smallskip

\noindent{{\bf {\em Case 2}.} {\em $z\in\OM$ is such that $\lmink{n\times n}(W)=0$.}}

\noindent{Let $t<1$ be so small that
\[
t \ = \  \left|\frac{\zt_2-\zt_1}{1-\zetbar_1\zt_2}\right|.
\]
Note that the map
\[
\fee_t(\zt) \ = \ \left[\frac{\zt W_{j,k}}{t}\right]_{j,k\leq n} \ = \ 
\frac{\zt W}{t}
\]
satisfies $\fee(D)\subset\OM_n$ since $\zt W\in\OM_n$ and $\lmink{n\times n}(W)=0$.
Write $(\zt_2-\zt_1)/(1-\zetbar_1\zt_2)=\alpha(\zt_1,\zt_2)$. Define
\[
\psi(\zt) \ := \ e^{-i{\rm Arg}(\alpha(\zt_1,\zt_2))}\frac{\zt-\zt_1}{1-\zetbar_1\zt}.
\]
It is obvious that $f:=\fee_t\circ\psi$ has the desired properties. This establishes
part (2) of Theorem \ref{T:specUB}.}
\end{proof}
\end{custom}

\end{document}